\documentclass[amsart]{ajmlatex} 
\usepackage{epsfig}
\usepackage{amscd}
\usepackage{latexsym}
\usepackage{amssymb}
\input amssym.def
\input amssym.tex
\def\bbbone{{\mathchoice {1\mskip-4mu \text{l}} {1\mskip-4mu \text{l}}
{ 1\mskip-4.5mu \text{l}} { 1\mskip-5mu \text{l}}}}

\newcommand{\vol}{\operatorname{vol}}

\newsymbol\circlearrowleft 1309
\newsymbol\restriction 1316

\renewcommand{\Re}{\mathop{\rm Re}\nolimits}
\renewcommand{\Im}{\mathop{\rm Im}\nolimits}

\newtheorem{thm}{Theorem}

\newcommand{\text}{\textrm}
\def\eqref#1{(\ref{#1})}

\begin{document}  
\slugline{AJM}{2}{3}{609--618}{September}{1998}{006}
\setcounter{page}{609}

\title{Poisson formula for resonances in even dimensions\thanks{Received
July 22, 1998; accepted for publication September 3, 1998.}}
\author{Maciej Zworski\thanks{Department of Mathematics,
University of Toronto, Toronto, ON M5S 3G3,
Canada (zworski@ math.toronto.edu) and Department of Mathematics,
University of California, Berkeley, CA 94720, USA 
(zworski@math.berkeley.edu).}}

\pagestyle{myheadings}
\thispagestyle{plain}
\markboth{MACIEJ ZWORSKI}{POISSON FORMULA FOR RESONANCES}
  
\maketitle  
  
\section{Introduction}  
We consider scattering by an abstract compactly supported perturbation    
in $ {\Bbb R}^n $. To include the traditional cases of potential,      
obstacle and metric scattering without going into their particular    
nature we adopt the ``black box" formalism developed     
jointly with Sj\"ostrand \cite{SjZw1}. It is quite likely that    
one could extend the results presented here to the case of     
non-compactly supported perturbation as well -- see \cite{Sj-trace}    
for a natural generalization of ``black box" perturbations.    
   
We review now the basic assumptions. We work with a complex Hilbert   
space with an orthogonal decomposition      
\begin{equation}   
\label{eq:bb1}   
{\cal H} = {\cal H}_{R_0 } \oplus L^2 ( {\Bbb R}^n \setminus B ( 0 ,    
R_0 ) ) \, , 
\end{equation}   
and with an operator    
\begin{eqnarray}   \label{eq:bb2}   
&& P \; : \; {\cal H} \longrightarrow {\cal H} \, , \ \    
\text{self-adjoint with a domain} \ {\cal D} \subset {\cal H}\nonumber \\   
&& {\bf 1}_{{\Bbb R}^n \setminus B ( 0 , R_0 ) } {\cal D} =    
H^2 ( {\Bbb R}^n \setminus B( 0 , R_0 ) ) \\   
&& {\bf 1}_{ {\Bbb R}^n \setminus B ( 0 , R_ 0 ) } P = - \Delta    
|_{ {\Bbb R}^n \setminus B ( 0 , R_0 ) }\,,   \nonumber
\end{eqnarray}   
which satisfies   
\begin{equation}   
\label{eq:bb3}   
\exists \; k \ \text{such that} \  \ {\bf 1}_{B ( 0 , R_0 )} ( P + i )^{-k}   
\ \text{ is of trace class }\,,    
\end{equation}   
\begin{equation}   
\label{eq:bb4}   
P \geq - C \, ,  \ \ \ C \geq 0 \,.   
\end{equation}   
These assumptions guarantee that  the resolvent $ R ( \lambda )   
= ( P - \lambda^2 )^{-1} $ continues meromorphically as an operator  
$ {\cal H}_{\rm{comp}} \rightarrow {\cal D}_{\rm{loc}} $ from $   
\Im \lambda < 0 $, $ \lambda^2 \notin \sigma_{\rm{pp}} ( P ) $, to   
$ {\Bbb C}  $ when   $ n $ is odd and to $ \Lambda $, the logarithmic  
plane, when $ n $ is even. The poles of this meromorphic continuation   
are called {\em resonances}. At $ \lambda \neq 0 $ all reasonable  
definitions of multiplicity agree. We can for instance say that  
the multiplicity of a pole at $ \lambda \neq 0 $ is given by   
the rank of the polar part of $ R ( \zeta ) $ near $ \lambda $ -- see  
\cite{Mel4}. The situation is more subtle at $ 0 $ and rather   
than go into a detailed discussion we will take the multiplicity   
required by the trace formula -- see \cite{Vai},\cite{Vod2},\cite{ZwP}
for the discussion of the resolvent near $ 0$. 
   
If $ U ( t) $ is the wave group for the operator $ P$ and $ U_0 ( t ) $    
is the free wave group we consider the natural {\em wave trace}:    
\[ u  ( t ) = \text{tr} \,  U ( t) - U_0 ( t)  \,, \] 
which is an even distribution in $ t \in {\Bbb R} $.  The notation 
used here is somewhat informal
since $ U $ and $ U_0 $ act on different 
spaces -- see \cite{SjZw5}. The correct definition is given by 
\begin{equation} 
\label{eq:u} \quad
u ( t )  \stackrel{\text{def}}{=} \text{tr}
\; \left( U(t) - 
\bbbone_{ {\Bbb R}^n \setminus   B(0, R_0 )}U_0 ( t )
\bbbone_{ {\Bbb R}^n \setminus B ( 0 ,R_0 ) }  \right)  
+ \text{tr}
\;
\bbbone_{  B(0, R_0 )}U_0 ( t )
\bbbone_{  B ( 0 ,R_0 ) }   \,. 
\end{equation}
In odd    
dimensions the following Poisson formula was established in     
increasing degrees of generality by Bardos-Guillot-Ralston \cite{BGR},    
Melrose \cite{Mel1},\cite{Mel2} and Sj\"ostrand-Zworski \cite{SjZw5}:    
\begin{equation}    \label{eq:otrace}    
t^{n+1} u  (  t )  =   t^{n+1}    
\sum_{\lambda \in {\Bbb C}}  m( \lambda )   e^{  i \lambda |t| } \,,     
\end{equation}

\vspace{-5mm}
$$
m( \lambda )  = \text{ multiplicity of $ \lambda $ as a resonance of $ P $}\,,  
$$
in the sense of distributions on $ {\Bbb R}$.     
The observation that we only need to multiply by $ t^{n+1} $ was made in     
\cite{ZwP}. The formula also holds exactly for super-exponentially     
decaying perturbations as was   
pointed out by S\'a Barreto-Zworski \cite{SaZw1}.      
    
We note that for $ t > 0 $ the trace formula is equivalent to     
\begin{equation}    
\widehat{ u \phi } ( \lambda ) = \sum_{\zeta\in {\Bbb C}} m ( \zeta )     
\hat \phi ( \lambda - \zeta)  \,, \ \ \phi \in {\cal C}^\infty_{\rm{c}}   
( ( 0 , \infty ) )  \,.    
\label{eq:otrace'}    
\end{equation}    
    
The original proofs of \eqref{eq:otrace} were based on Lax-Phillips    
theory \cite{LP} and in particular on the strong Huyghens principle.     
The extension to the case of hyperbolic surfaces by Guillop\'e-Zworski    
\cite{GuZw3} provided a proof which does not require the strong Huyghens     
principle and is also applicable in the euclidean case \cite{ZwP}. 
It is based on the Birman-Krein formula and ``global minimum modulus"    
estimates on the scattering determinant. That was followed by a local     
trace formula of Sj\"ostrand \cite{Sj-trace} the proof of which did not    
involve any scattering theory but also used some ``local minimum    
modulus" estimates for determinants of some holomorphic matrices.     
Sj\"ostrand's formula specialized to  
the even dimensional compactly supported    
case gives the following weaker version of \eqref{eq:otrace'}:    
\begin{equation}    \label{eq:etrace'}    
\widehat{ u \phi } ( \lambda )  = 
\sum_{\zeta\in \lambda \Omega} m ( \zeta ) \hat \phi ( \lambda - \zeta)      
+ {\cal O } ( \langle \lambda \rangle^{-\infty } )     \,,
\end{equation}

\vspace{-5mm}
\[
\Omega = [ 1/2, 3/2] + i [0, 1/2]\,, \ \     
   \phi \in {\cal C}^\infty_{\rm{c}}((0, \infty )) \,.      
\]
We remark however that the semi-classical local formula of \cite{Sj-trace}
is much stronger 
than \eqref{eq:etrace'}. 
    
By using a ``local minimum modulus" theorem   
in the     
argument of \cite{GuZw3}, \cite{ZwP} we can strengthen \eqref{eq:etrace'}    
to obtain a global formula. This extension was motivated    
by a question asked by Vodev (see Sect.3).  
  
\begin{thm}  
Let $ P $ be an operator sastisfying the     
assumptions \eqref{eq:bb1}-\eqref{eq:bb4} and let $ u  (  t ) $     
be its normalized wave trace given by \eqref{eq:u}.
Let $ \Lambda_\rho $ be an open conic    
neighbourhood of the real axis as shown in Fig.1, $ \sigma (     
\lambda ) $ the scattering phase of $ P $ and let  $ \psi \in     
{\cal C}^\infty_{\rm{c}} ( {\Bbb R} ; [0,1] ) $ be equal to     
$ 1 $ near $ 0 $. Then     
\begin{eqnarray}    \label{eq:etrace}    
u (t) & = & \sum_{ \lambda \in \Lambda_\rho } m ( \lambda ) e^{i \lambda |t| } 
 + \sum_{ {\lambda^2 \in \sigma_{\rm{pp}} (P)  \cap (-\infty, 0 )  
 } \atop {  
\Im \lambda < 0 } } m ( \lambda ) e^{ i  \lambda |t|}   \\
 & &  + \  m( 0 ) +   
2 \int_{ 0}^\infty \psi ( \lambda ) \frac{d \sigma}{d \lambda}     
( \lambda ) {\cos  t  \lambda } d \lambda + v_{ \rho , \psi } ( t) \,,
\ \  t \neq 0 \,, \nonumber \\ 
v_{\rho, \psi } & \in & {\cal C}^\infty ( {\Bbb R} \setminus \{0 \}) \,,    
\ \   \partial^k_t  v_{\rho, \psi} = {\cal O} ( t^{-N} )\,, \ \forall \; 
k, N \,, \ \      
|t| \longrightarrow \infty \,.     \nonumber
\end{eqnarray}    
\end{thm}  
  
\begin{figure}[htb]  
\centerline{\psfig{figure=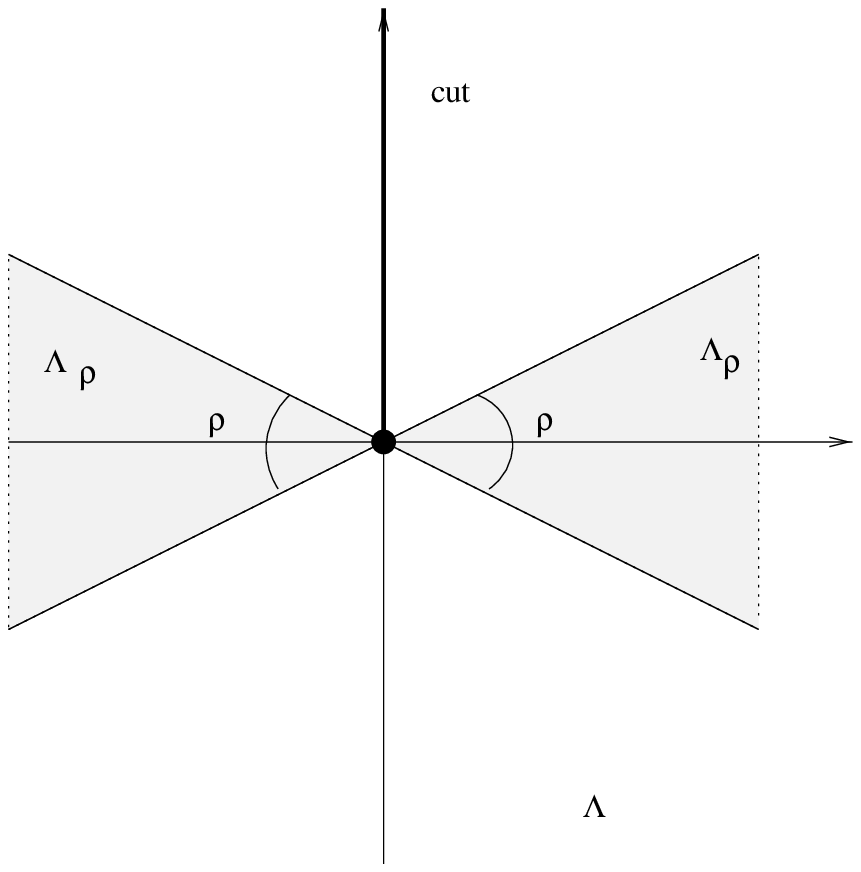,height=9cm}}
\centerline{{\small {\sc Fig. 1.} 
{\it Conic neighbourhoods of the real axis on the logarithmic  
plane.} } } \label{Fig.1}  
\end{figure}  
  
The scattering phase, $ \sigma ( \lambda ) $, is a standard object  
in scattering theory -- see \cite{Mel4}  and reference given there
for background information and \cite{Chr} for the discussion of the
``black box'' case.  
Here it is normalized so that the Birman-Krein formula holds -- see   
\eqref{eq:bikr} below.  
  
\section{Proof of the trace formula}  
To prove Theorem 1 we identify the subset of $ \Lambda $    
shown in Fig.1 with $ {\Bbb C} \setminus e^{ i \pi /2 } 
\overline{\Bbb R}_+ $,   
where $ e^{  i \pi / 2 } {\Bbb R} $ is the cut. Then the resonances    
are symmetric with respect to the cut and they coincide with the    
poles of the scattering determinant $ s ( \lambda ) = \det    
S( \lambda ) $ . The unitarity of $ S ( \lambda ) $ for $ \lambda > 0   
$  implies the usual relation $ S ( \lambda )^{-1}  =  S(    
\bar \lambda )^* $ for $ \Re \lambda > 0 $.   
Hence $ s ( \lambda )^{-1} = \overline{s (    
\bar \lambda )} $. To simplify the discussion\footnote{so that  
we do not need to consider global analytic properties of $ S (   
\lambda ) $} we will consider  
the scattering matrix in $ \Lambda_\rho \cap \{ \Re \lambda > 0 \}  
$ only and {\em define} $ s ( \lambda ) $ in $ \Lambda_\rho   
\cap {\Re \lambda <  0 } $ so that the scattering phase defined   
by $ \sigma' ( \lambda ) = ( i/ 2 \pi ) s' ( \lambda ) / s ( \lambda ) $  
is even in $ \lambda $ -- see \eqref{eq:bikr}.  
  
The assumption \eqref{eq:bb3} guarantees the existence of $ m $ such   
that for  $ \rho > 0 $  
\begin{equation}   
\label{eq:vod}   
\sum_{{\lambda \in \Lambda_\rho} \atop {|\lambda | \leq r }} 
m( \lambda ) \leq C_\rho r^{m  
+ \epsilon}  \,, \ \ \forall \; \epsilon > 0 \,,  
\end{equation}   
see \cite{ZwP}. This is deduced from   
the polynomial bounds of Vodev \cite{Vod2},\cite{Vod3}    
which for $ \rho < \pi /2 $ (all that is needed here) follow also    
from the earlier estimates of Sj\"ostrand-Zworski \cite{SjZw1}.   
  
We now put  
\begin{equation}   
\label{eq:Weier}    
P_{\rho } ( \lambda ) = \prod_{\zeta \in {\Lambda_{\rho}} \setminus  
{\Bbb R} }   
E \left(  {\lambda}/{\zeta } , m  \right)^{m( \zeta )} \,, \ \    
E ( z , p ) = ( 1 - z ) e^{ z + \frac{z^2}2 + \cdots + \frac{z^p}{p}}    
\,,    
\end{equation}   
where the bound \eqref{eq:vod} guarantees the convergence of the Weierstrass   
product.  
   
The poles of $ s ( \lambda ) $ and $ s (  \lambda )^{-1} $ in $ \Lambda_\rho $   
coincide (with multiplicities) with the zeros of $ P_\rho ( \lambda ) $    
and $ P_\rho  ( - \lambda ) $ respectively.  Hence we can write   
\begin{equation}   
\label{eq:fact}   
s( \lambda ) = e^{ g _\rho ( \lambda ) } \frac{P_\rho ( - \lambda ) }{   
P_\rho  (  \lambda ) } \,, \ \ \  \lambda \in \Lambda_\rho \cap   
\{ \Re \lambda > 0 \} \,,   
\end{equation}   
where $ g_\rho $ is holomorphic in $ \Lambda_\rho \cap   
\{ \Re \lambda > 0 \}$. We now extend $ g _ \rho ( \lambda ) $, and  
consequently $ s ( \lambda ) $,  to $ \Lambda_\rho $ by setting  
\[ g_\rho ( - \lambda ) = - g _\rho ( \lambda ) \,.\]  
That clearly implies that $ s' ( \lambda ) / s ( \lambda ) =   
s' ( - \lambda ) / s (- \lambda ) $ for $ \lambda \in {\Bbb R }  
\setminus \{ 0 \} $ and further analysis shows that this identity   
holds through $ \lambda = 0 $ -- see Sect.3.  
  
We want to estimate the function $  g_\rho $. For that we need   
to estimate $ s ( \lambda ) $ away from its poles and that is done   
exactly as in \cite{GuZw3},\cite{ZwP} (see also \cite{SaTa2}). We write   
\begin{eqnarray}  \label{eq:A} 
& & S ( \lambda )  = Id + A ( \lambda ) \,, \\    
& & A ( \lambda ) = C_n  \lambda^{n-2}    
{\Bbb E}^{\phi_1 } ( - \lambda ) ( I + K ( \lambda , \lambda_0 ) )^{-1}   
[ \Delta, \chi ] {}^t {\Bbb E}^{\phi_2 } ( \lambda ) \,, \nonumber  \\  
& & {\Bbb E}^\rho : L^ 2 ( {\Bbb R}^n ) \longrightarrow L^2 ( {\Bbb S}^{n-1} )\,,  
\ \ {\Bbb E}^\rho ( \theta, x ) = e^{ i \lambda \langle x , \theta \rangle }  
\rho ( x ) \,, \ \rho \in {\cal C}^\infty_{\rm{c}} ( {\Bbb R}^n ) \,, \nonumber
\end{eqnarray}   
and where $ K ( \lambda , \lambda_0 ) $ is the operator constructed in    
Sect.3 of \cite{SjZw1}:   
\begin{eqnarray}   \label{eq:K} 
\qquad K( \lambda , \lambda_0  ) = [\Delta , \chi_0 ] R_0 ( \lambda )  ( 1 -    
\chi_1 ) \chi_4   
- [\Delta , \chi_2] R ( \lambda_0 ) \chi_4  + \chi_2 ( \lambda_0 ^2    
- \lambda^2 ) R ( \lambda_0 ) \chi_4 \,, \\   
\chi_i \in {\cal C}_{\rm{c}}^\infty ( {\Bbb R}^n ) \,, \    
\chi_0 \equiv 1 \ \text{near $ B ( 0 , R_0 ) $}\,, \    
\chi_i \equiv 1 \ \text{near $ \text{supp}\; \chi_{i-1} $}\,, \nonumber \\   
R ( \lambda_0 ) = ( P - \lambda_0^2 )^{-1} \,, \ \Im \lambda_0 \ll 0 \,, \ \   
R_0 ( \lambda ) = ( - \Delta - \lambda^2   )^{-1} \,.    \nonumber
\end{eqnarray}   
To estimate $ s ( \lambda ) $ we will first estimate $ \|    
( I + K ( \lambda, \lambda_0 ) )^{-1} \|$ and that is based on    
the inequality    
\begin{equation}   
\label{eq:gokr}   
  \| ( I + K ( \lambda, \lambda_0 )) ^{-1} \| \leq    
\frac{ \det (  I + | K ( \lambda , \lambda_0 )|^{m+1} )  } {   
| \det (I + K ( \lambda , \lambda_0) )^{m+1 } ) |}  \,,    
\end{equation}   
from \cite{GoKr}, Theorem 5.1, Chap.V.   
   
Exactly as in \cite{ZwP}, where we followed \cite{Vod1},\cite{Zw1},   
we see that for $ \lambda \in \Lambda_\rho $    
\[ | \det  ( I + K ( \lambda , \lambda_0 )  ^{m+1}  ) | \leq    
C e^{ C|\lambda|^{m + \epsilon}}  \,. \]   
Using the lower modulus theorem of H. Cartan\footnote{We remark    
that a much cruder estimate would suffice here but it is nice   
to quote the optimal result which is useful elsewhere   
in the theory of resonances -- see \cite{SaTa2} and 
\cite{SjZw7}, Sect.8.} 
 -- see for instance \cite{levin}, Theorem 4, Sect.11.3 --   
we obtain a lower bound:   
\[ | \det ( I +  K ( \lambda , \lambda_0 )  ^{m+1}  ) | \geq    
C   e^{ - C r ^{m + \epsilon}/ \eta }  \,, \ \ \ \lambda \in    
D( r , \rho r /C ) \setminus  \bigcup_{j} D(    
\lambda_j , r_j ) \,, \ \  \sum_j r_j \leq \eta r  \,, \]   
uniformly as $ r \rightarrow \infty  $. From this and    
\eqref{eq:gokr} we obtain, as in \cite{ZwP},   
\[ | s ( \lambda ) | \leq C e^{ C  (r^{m + \epsilon}/\eta )^n }    
\,, \ \ \ \lambda \in    
D( r , \rho r /C ) \setminus  \bigcup_{j} D(    
\lambda_j , r_j ) \,, \  \ \sum_j r_j \leq \eta r  \,. \]   
If we take $ \eta \ll 1  $ then for every    
$ r $ there exists $ r/2C \leq k(r) \leq  r/C $ such that the  
circle   
$ | \lambda - r| = k ( r) $ does not intersect any of the    
excluded discs. Then using the standard estimates for    
Weierstrass products and the maximum principle we see that   
\[ | \exp { g_\rho ( \lambda ) } | | P_\rho ( - \lambda ) |    
\leq C e^{ |\lambda|^{ (m +\epsilon)
n } } \,, \  \lambda \in    
\Lambda_{\rho'} \,, \ 0 < \rho' \ll \rho \,. \]   
We then conclude (as in the proof of Cartan's theorem or    
yet easier as in the proof of Hadamard's factorization    
theorem\footnote{We can use for instance Carath\'eodory's    
inequality -- see \cite{Tits}.}) that    
\begin{equation}   
\label{eq:sym} |\partial^k_\lambda g_{\rho}{ ( \lambda ) } | \leq C |   
\lambda |^{ (m  + \epsilon )
n - k } \,, \ \ \lambda \in    
\Lambda_{\rho''} \,, \ \ 0 < \rho'' < \rho' \,, \end{equation}   
where the symbolic property followed from Cauchy's inequalities.   
   
We can now prove the Poisson formula. As in \cite{GuZw3},\cite{ZwP}   
the starting point is the Birman-Krein formula:   
\begin{equation}   \label{eq:bikr}   
u (t) = \widehat{ \sigma' }  ( t)  +    
\sum_{ \lambda^2 \in \sigma_{\rm{pp}} ( P ) \setminus \{ 0 \} }  
2 \cos ( t \lambda ) + m ( 0 ) \,, 
\end{equation}

\vspace{-6mm}
\[
\sigma ( \lambda ) \stackrel{\text{def}}{=}  
 \frac{i}{2\pi}  \log s ( \lambda )   \   
\text{for} \ \lambda > 0 \,, \ \ \sigma ( \lambda ) =  - \sigma ( - \lambda )   
\ \text{for} \ \lambda < 0 \,,  
\]
where the Fourier transform is of course taken in the sense   
of distributions. For the ``black box" perturbation the proof was    
given by Christiansen in Sect.1 of \cite{Chr}
but it is classical for all well known scattering problems. We note  
that our definition of $ s ( \lambda ) $ implies that we can set  
\[ \sigma' ( \lambda ) = \frac{i} { 2 \pi } \frac{ s' ( \lambda ) }  
{s ( \lambda ) } \,, \ \ \forall \; \lambda \in {\Bbb R}   
\setminus \{ 0 \}\,. \]  
   
Let $ \psi $ be as in Theorem 1 and even. We can then write   
\begin{eqnarray*} 
\sigma' ( \lambda )  & = &  \frac{i}{2\pi}{\frac{ s' (\lambda)}{s(\lambda)}}\\
& = &    
\frac{ i }{ 2 \pi } \left( (1 -    
\psi ( \lambda )) g'_\rho ( \lambda ) + \frac{d}{d \lambda}   
\left( \log P_\rho (  - \lambda )   
-  \log P_\rho ( \lambda ) \right) \right. \\
& & \qquad \qquad \left.
 + \psi ( \lambda ) \frac{s' ( \lambda ) }{   
s ( \lambda ) }  \right) +  f_{\rho , \psi} ( \lambda ) \,, 
\end{eqnarray*}
where $ f_{\rho, \psi } \in {\cal C}_{\rm{c}}^\infty ( {\Bbb R} ) $.  
  
The argument of \cite{GuZw3},\cite{ZwP} now easily gives   
\begin{eqnarray*}
t^{2p } u  (  t ) & = & t^{2p} \left( \sum_{  
\lambda \in \Lambda_\rho } e^{ i \lambda |t| } m   
( \lambda )   \right . \\
& & \quad  \left.
+ \sum_{ {\lambda^2 \in \sigma_{\rm{pp}} (P)  \cap (-\infty, 0 )  
 } \atop {  
\Im \lambda < 0 } } e^{ i  \lambda |t|}  + m( 0 ) +  \widehat{ \psi 
\sigma' } ( t)   + 
\widehat{ ( 1 - \psi) g_\rho } ( t)  + \hat{f}_{{\rho, \psi}} (t) \right)  
\,, 
\end{eqnarray*}   
for $ p $ large enough.  
Since $ g_\rho $ is a symbol on $ {\Bbb R} $,  
\[    v_{\rho, \psi} ( t ) =   
\frac{i}{ 2 \pi }   
\widehat{ ( 1 - \psi) g_\rho } ( t)  + \hat{f}_{\rho, \psi} (t)  \,, \]   
has the properties stated in Theorem 1 and this completes its proof.  
   
We remark here that {\em a posteriori} the bound on $ g_\rho $ on the   
real axis has to be much better than the bound provided by the estimate  
\eqref{eq:sym}: we know the strength of the singularity of $ u $ at   
$ t = 0 $ and the bound on the number of resonances gives an estimate  
on the strength of the singularity of the exponential sum. Hence for  
elliptic perturbations where $ m = n $ we obtain  
\[ |\partial^k_\lambda  g_\rho ( \lambda ) | \leq C_{k ,   
\epsilon} ( 1 + |\lambda| )^{n + \epsilon  - k }  \,, \ \   
\forall \, \epsilon > 0 \,. \]  
      
\section{Review of applications}   
The basic application of the trace formula is in obtaining  
lower bounds on the number of resonances from the singularities  
of the wave trace. The basic Tauberian lemma  
was given in Sj\"ostrand-Zworski \cite{SjZw3} and it was applied there  
to problems in odd dimensions. One of the applications of the  
local trace formula of Sj\"ostrand \cite{Sj-trace} was the extension
of those results to even dimension -- see Theorem 10.1 there.   
That becomes even clearer when we use the global formula  
\eqref{eq:etrace}. One of the interesting consequencies is based on the   
trace formula of Guillemin-Melrose \cite{GuMe}:  
\begin{thm}  
Let $ P $ be the Dirichlet or Neumann Laplacian on a connected exterior   
domain $ {\Bbb R}^n   
\setminus {\cal O} $ where $ {\cal O } $ has a smooth boundary.  
Suppose that the there exists a non-degenerate closed transversally  
reflected   
trajectory of the broken geodesic flow of $ {\Bbb R}^n   
\setminus {\cal O} $ such that no essentially different closed trajectory   
has the same period. Then for any $ \epsilon > 0 $ there exists $ C_\epsilon 
> 0 $ such that   
\[ \sum \{ m( \lambda ) \; :\; |\lambda | \leq r + C_\epsilon \,, \   
|\Im \lambda | < \epsilon  \log |\lambda | \} \geq r / C_\epsilon \,. \]  
\end{thm}  
\noindent 
For more applications we refer to \cite{SjZw3} and \cite{Sj-trace},  
Sect.10.  
  
The new trace formula allows also an easy extension of   
some of   
the results of Ikawa on the distribution of resonances for several  
convex obstacles to even dimensions. That is particularly   
interesting in dimension two where most of the numerical studies  
were conducted -- see for instance \cite{num} for the discussion 
of symbolic dynamics. We remark that it is rather clear  
that the results of \cite{Sj-trace} would suffice for this   
purpose but the global formula makes the applications even   
more apparent. As an example we give the modification of the  
result of \cite{Ik}:  
\begin{thm}  
Let $ P $ be the Neumann Laplacian on $ {\Bbb R}^n   
\setminus {\cal O} $, $ n \geq 2 $, and let us assume that $  
{\cal O} = \bigcup_{j=1}^N {\cal O}_j $ where $ {\cal O}_j $  
are mutually disjoint strictly convex obstacles with smooth boundaries  
satisfying the following   
condition:  
\[  \text{\em convex hull } ( {\cal O}_k \cup {\cal O}_l )  
\cap {\cal O}_m = \emptyset \ \ \ \forall \; k \neq m \neq l\,.\]  
Then there exists $ \alpha > 0 $ for which  
\[ \sum \{ m( \lambda ) \; : \; |\Im \lambda| < \alpha \} =  \infty\,.  
\]  
\end{thm}  
  
The next theorem  answers a question asked by Vodev \cite{VodP}  
and does not seem to follow from the local trace formula:  
  
\begin{thm}  
Let $ P $ satisfy the assumptions \eqref{eq:bb1}-\eqref{eq:bb3}  
with $ n \geq 4 $ even.  
In addition let us assume that $ 0$   
is not an eigenvalue or a zero-resonance  
of $ P$. Then for any $ \gamma > 0 $  and $ k$  
\begin{equation}  
\label{eq:est}   \qquad
\left|\left(\frac{ \partial} {\partial t} \right)^{ k } 
 \left( u  (  t ) - \sum_{ \Im \lambda  \leq \gamma   
\log |\lambda | } m( \lambda ) e^{ i \lambda |t| } \right) \right| \leq
C_k  t^{ -n+2 -k } \,, \ \ t > t_k >   (n+k) / \gamma \,. \end{equation}  
\end{thm}  
We remark that in odd dimensions the polynomial bound on the  
number of resonances and the global formula \eqref{eq:otrace} 
imply that the right hand side of \eqref{eq:est} can be  
replaced by $ {\cal O} ( \exp ( - \alpha t ) )$, $ \alpha > 0 $. 
It is quite clear from \eqref{eq:etrace} that to establish Theorem 4 
we need to understand the behaviour of $ \sigma' ( \lambda ) $ as 
$ \lambda \rightarrow 0 + $. A finer analysis based for instance on 
\cite{J},\cite{JK}, should show that the estimate above is optimal  
and that in fact there exists an asymptotic expansion
as $  t \rightarrow \infty $.  
Also, we did not attempt to study the 
more  
involved two dimensional case. 
\begin{proof}  
We start by observing that  
\begin{eqnarray*}
\left| \left(\frac{ \partial} {\partial t} \right)^{ k } 
 \left(  \sum_{ \Im \lambda  \leq  \rho |\lambda | }  
m( \lambda ) e^{ i \lambda |t| } 
-  \sum_{ \Im \lambda  \leq \gamma   
\log |\lambda | } m( \lambda ) e^{ i \lambda |t| } \right)   \right|
\leq C_k  e^{ - \alpha_k t }  \,, \\ \alpha_k > 0 \,, 
\ \ t \geq t_k > (n+k) /\gamma \,, 
\end{eqnarray*} 
which follows from the bound \eqref{eq:vod}. In fact, the derivative
of the difference can be estimated by 
\begin{eqnarray*}
 e^{ - \alpha'_k t } + \int_A ^\infty x^k x^{ - \gamma t } d N( x) 
& \leq& e^{ - \alpha' _k t } + C A^{ - \gamma ( t - (n+k)/\gamma )}
+ C t \int_A^\infty x^{ - \gamma ( t - (n+k)/\gamma ) -1  } dx 
\\ & \leq& C_{k, 
\epsilon} e^{ - \alpha_k t } \,, \ \ 
 t  > ( n + k ) / \gamma + \epsilon \,. 
\end{eqnarray*}
 
Thus we need to study the behaviour of  
\begin{equation} 
\label{eq:strace} 
\frac {s' ( \lambda ) }{ s ( \lambda ) } = \text{tr} \; S( \lambda   )^* 
S' ( \lambda ) = \text{tr} \; ( I + A ( \lambda )^* ) A' ( \lambda ) \,, 
\end{equation} 
as $ \lambda \rightarrow 0+ $, where we used the notation of \eqref{eq:A}. 
 
The assumption that $ P $ has no resonance at zero implies that the 
cut-off resolvent, $ \chi R ( \lambda ) \chi$, $ \chi \in {\cal C}_{\rm{c}}^ 
\infty ( {\Bbb R}^n )$, $ \chi \equiv 1 $ near $ B ( 0 , R_0 ) $,  
is  holomorphic in $ ( \lambda , \lambda^{n-2} \log \lambda ) $ for 
$ |\lambda | < \epsilon $.  To see that we recall that the free 
resolvent, $ R_0 ( \lambda ) $, is of the form $ R'_0 ( \lambda ) +  
\lambda^{n-2} \log \lambda M ( \lambda ) $ where $ R' _0 ( \lambda ) $ 
and $ M( \lambda ) $ are entire -- see Sect.1 of \cite{Mel4}.  
Following \cite{SjZw1}, Sect.3, we write $ \chi R ( \lambda ) \chi 
= \chi ( Q_0 ( \lambda ) \chi + Q_1 ( \lambda_0 ) \chi ) ( I + K  
( \lambda , \lambda_0 ))^{-1} $ where $ Q_0 $ is a cut-off free  
resolvent and $ K $ is given by \eqref{eq:K}. Analytic Fredholm 
theory shows that when $ \chi R ( \lambda ) \chi $ is bounded, it is 
a holomorphic function of $ \lambda $ and $ \lambda^{n-2} \log \lambda$  
-- see \cite{Vai}, \cite{Vod2} for more details.  
 
To study \eqref{eq:strace} we need a different representation of  
$ A ( \lambda ) $. If we recall the definition from Sect.3 of \cite{ZwP}, 
$ A ( \lambda ) $ comes from the radiation pattern of $ R (  
\lambda ) ( -  [ \Delta , \chi ] e^{ i \langle \bullet, \omega \rangle }  
) $. To obtain a formula similar to \eqref{eq:A} but involving $ R (  
\lambda ) $ rather than $ ( I + K ( \lambda , \lambda_0 ) ) ^ { -1} $ 
we take $ \chi_j $ as in \eqref{eq:K} and write  
\[ (  1 - \chi_2) R ( \lambda ) \chi_1 =  
R_0 ( \lambda ) (-  \Delta - \lambda^2 ) ( 1 - \chi_2 ) R ( \lambda )  
\chi_1 = R_0 ( \lambda ) ( - [ \Delta , \chi_2 ] R ( \lambda ) \chi_1 )\,,\] 
since $  ( 1 - \chi_2 )  ( - \Delta - \lambda ^2 ) 
=  ( 1 - \chi_2 )  ( P - \lambda^2 ) $ and $ ( 1 - \chi_2 ) \chi_1 = 0 $. 
This shows that 
\[ A ( \lambda )  = c_n \lambda^{n-2}  
{\Bbb E}^{\phi_1 } ( - \lambda ) [ \Delta, \chi_2 ] R ( \lambda )  
[ \Delta, \chi ] {}^t {\Bbb E}^{\phi_2 } ( \lambda ) \,,   \\  \]
and the assumption on $ R ( \lambda ) $ implies that $  
\lambda ^{ 2 - n } A ( \lambda ) $ is holomorphic in $ \lambda $ and 
$ \lambda^{n-2} \log \lambda $ for $ |\lambda| < \epsilon $.  
 From this and \eqref{eq:strace} it follows that  
\begin{equation} 
\label{eq:f} 
 \sigma'  ( \lambda ) = \lambda^{n-3} f ( \lambda , \lambda^{n-2} 
\log \lambda ) \,,  \ \ \lambda > 0 \,, \end{equation} 
with $ f $ smooth near $ 0 $. We then easily check that 
\begin{eqnarray*}
& & \int_0^\infty \lambda^{n-3} f ( \lambda , \lambda^{n-2}  
\log \lambda ) \psi ( \lambda ) \cos t \lambda d \lambda   \\
& = &(-1)^{\frac{n}2} ( n-3) ! f ( 0,0) \; t^{-n+2} + {\cal O} ( t ^{ 
-n+1} ) \,,  \ \ t \rightarrow \infty \,, 
\end{eqnarray*}
which completes the proof as the estimates for derivatives clearly  
hold as well. 
\end{proof} 
 
Finally we compare this result with the estimates on the heat  
trace. As a consequence of well known estimates on heat kernels, 
S\'a Barreto-Zworski \cite{SaZw2} showed that when $ P =  
- \Delta + V $ and $ P $ has no resonances with $ \Im \lambda  
\leq  0 $ (that is no eigenvalues and no zero resonance) 
\[ \text{tr} \; ( e^{ - t P } - e^{ t \Delta } ) = {\cal O} (  
t ^{ -\frac{n}2 + 1 }) \,, \ \ t > 0 \,. \] 
Werner M\"uller pointed out to the author that for the behaviour 
as $ t \rightarrow \infty $ it is more natural to study  
$ \sigma' ( \lambda ) $ near $ \lambda = 0 $ using 
the heat version of the Birman-Krein formula: 
\begin{equation} 
\label{eq:birk'} 
\text{tr} \; ( e^{ - t P } - 
  e^{ t \Delta }
  ) =  
\int_0^\infty e^{ - t \lambda^2 } \sigma'( \lambda ) d \lambda 
+ \sum_{ \mu_j \in \sigma_{\rm{pp}} ( P ) } e^{ - t \mu_j }  
\,,  
\ \  t > 0 \,, 
\end{equation}
where to make sense of the trace 
we used the convention employed in the definition of 
$ u $, \eqref{eq:u}.
Hence under the assumptions of Theorem 4 but for {\em any} $ n \geq 3$ 
we obtain from its proof 
\begin{equation} 
\label{eq:sazw} 
\text{tr} \; ( e^{ - t P } - e^{ t \Delta } ) =  
 \sum_{{ \mu_j \in \sigma_{\rm{pp}} ( P ) } 
\atop{ \mu_j < 0 } } e^{ - t \mu_j }  +  
\frac12 \Gamma \left( \frac{n}2 -1  \right) f ( 0 , 0 ) \; t^{ -\frac{n}2  
+ 1 } + {\cal O} ( t^{ -\frac{n}{2} + \frac12} ) \,, 
\end{equation} 
where $ f $ is as in \eqref{eq:f}. In odd dimension it is a function 
of one variable, $ \lambda $, only.

\vspace{2mm}
{\sc Acknowledgments.} I should like to thank Laurent Guillop\'e and
Georgi Vodev for helpful 
comments on the first version of this paper. The  
partial support of this work  by the National Science
and Engineering Research Council of Canada and by the Erwin
Schr\"odinger Institute is also gratefully acknowledged.

\newpage
\mbox{ }
  
\end{document}